\documentclass[a4paper,12pt]{amsart}
\usepackage[utf8]{inputenc}
\usepackage[T1]{fontenc}
\usepackage[UKenglish]{babel}
\usepackage[a4paper,margin=28mm]{geometry}
\usepackage{verbatim}

\allowdisplaybreaks[4]
\usepackage{times}
\usepackage{dsfont,mathrsfs}
\usepackage{amsmath}
\usepackage{amsthm}
\usepackage{amssymb}
\usepackage{amsfonts}
\usepackage{latexsym}
\usepackage{booktabs}

\usepackage[colorlinks,linkcolor=red,pagebackref]{hyperref}
\usepackage{xcolor}

\usepackage{ulem}

\newtheorem{theorem}{Theorem}[section]
\newtheorem{lemma}[theorem]{Lemma}

\newtheorem{assumption}[theorem]{Assumption}

\theoremstyle{definition}
\newtheorem{definition}[theorem]{Definition}

\newtheorem{remark}[theorem]{Remark}
\numberwithin{equation}{section}
\renewcommand{\labelenumi}{\roman{enumi})}
\renewcommand\theenumi\labelenumi
\renewcommand{\leq}{\leqslant}

\renewcommand{\geq}{\geqslant}

\newcommand{\tl}{\tilde}
\newcommand{\Be}{\begin{equation}}
\newcommand{\Ees}{\end{equation*}}
\newcommand{\Bes}{\begin{equation*}}
\newcommand{\Ee}{\end{equation}}

\newcommand{\R}{\mathbb{R}}
\newcommand{\E}{\mathbb{E}}

\newcommand{\PP}{\mathbb{P}}

\newcommand{\Gcal}{\mathcal{G}}
\newcommand{\Lcal}{\mathcal{L}}

\newcommand{\dif}{\mathrm{d}}

\usepackage{marginnote}
\begin{document}
\title[Convergence rate of randomized midpoint Langevin Monte Carlo ]
{Convergence rate of randomized midpoint Langevin Monte Carlo}

\author[R. Li]{Ruinan Li}
\address{Ruinan Li: School of Statistics and Data Science, Shanghai University of International Business and Economics, Shanghai, China;}
\email{ruinanli@amss.ac.cn}

\author[T. Shen]{Tian Shen}
\address{Tian Shen: School of Statistics and Data Science, Shanghai University of International Business and Economics, Shanghai, China;}
\email{shentian@suibe.edu.cn}

\author[Z. Su]{Zhonggen Su*}
\address{Zhonggen Su: School of Mathematical Sciences, Zhejiang University, Hangzhou, Zhejiang, China;}
\email{suzhonggen@zju.edu.cn}


\begin{abstract}
The randomized midpoint Langevin Monte Carlo (RLMC),  introduced by Shen and Lee (2019),  is a variant of classical Unadjusted Langevin Algorithm. It was shown in the literature that the RLMC is an efficient algorithm for approximating high-dimensional probability distribution $\pi$. In this paper, we  establish the exponential ergodicity of RLMC with constant step-size.  Moreover, we design a dereasing-step size RLMC   and provide its convergence rate in terms of a functional class distance.

\end{abstract}

\maketitle


\section{Introduction} \label{Model}
Let $\pi$ be a  probability distribution on $\R^d$ with the form  $\pi(\dif x)=\frac{1}{Z}\exp(-U(x))\dif x$,
where $U$ is called the potential function and $Z=\int_{\R^d}\exp(-U(x))\dif x$ is the normalized  constant which  is unknown in general.
Sampling from such a target distribution $\pi$   is a fundamental algorithmic problem with applications in computational statistics, engineering, and machine learning; see for example the textbook Chewi \cite{che}. The standard method is Markov Chain Monte Carlo: construct a Markov chain whose stationary distribution is close to the target $\pi$ and run it until approximate stationarity. Consider the overdampted  Langevin diffusion
\begin{align}\label{ld}
\dif X_t=-\nabla U(X_t)dt+\sqrt{2}dB_t, \quad X_0\in\R^d,
\end{align}
where $(B_t)_{t\geq0}$ is a $d$-dimensional Brownian motion. Under mild regularity conditions on $U$, the solution $(X_t)_{t\geq0}$ of \eqref{ld} admits $\pi(\dif x)=\frac{1}{Z}\exp(-U(x))\dif x$ as its unique stationary distribution. In the literature, the common approach to construct a Markov chain for the approximation of $\pi$ involves the Euler-Maruyama discretization of Langevin diffusion \eqref{ld}: given step size $\eta$,   define
\begin{align*}
    \tl{X}_{(k+1)\eta}= \tl{X}_{k\eta}-\eta\nabla U( \tl{X}_{k\eta})+\sqrt{2\eta}\xi_{k+1}, \quad k\geq0,
\end{align*}
where $\{\xi_k\}_{k\geq1}$ is a sequence of $d$-dimensional standard Gaussian random vectors and . This recursion defines the famous Unadjusted Langevin Algorithm (ULA). Under suitable conditions, the ULA admits a unique stationary distribution $\mu_\eta$, which is close to the target $\pi$ when the step size $\eta$ is small enough. The convergence rate of the ULA to $\pi$ has been extensively studied in the literature, such as \cite{che, da, du, du2, fa, pages} and references therein.

 To enhance the algorithm's efficiency, Shen and Lee \cite{sh} proposed the randomized Langevin Monte Carlo (RLMC): for $k\geq0$
\begin{align}\label{rmm}
X_{(k+1)\eta}&=X_{k\eta}-\eta\nabla U(X_{(k+u_k)\eta})+\sqrt{2\eta}\xi_{k+1},
\end{align}
and
\begin{align*}
X_{(k+u_k)\eta}&=X_{k\eta}-u_k\eta\nabla U(X_{k\eta})+\sqrt{2u_k\eta}\xi_{k+1}^\prime,
\end{align*}
where $(u_k)_{k\in\mathbb N_0}$ is a sequence of i.i.d. uniform random variables on the interval $[0,1]$, $(\xi_k)_{k\geq1}$ and $(\xi_k^\prime)_{k\geq1}$ are sequences of i.i.d $d$-dimensional standard Gaussian vectors independent of $(u_k)_{k\in\mathbb N_0}$. According to \eqref{rmm}, $(X_{k\eta})_{k\geq0}$ is a Markov chain on the state space $\R^d$ with the transition density 
function
\begin{align}\label{zymd}
q_{\eta}(x,\tl x)=&\int_0^1\int_{\R^d}\frac{1}{(4\eta\pi)^{d/2}}\exp\left(-\frac{[\tl x-(x-\eta\nabla U(y))]^2}{4\eta}\right)\nonumber\\
&\qquad \cdot \frac{1}{(4\eta u\pi)^{d/2}}\exp\left(-\frac{[y-(x-u\eta\nabla U(x))]^2}{4u\eta}\right)\dif y\dif u
\end{align}
and the transition kernel 
\begin{align}\label{zyh}
Q_\eta(x,A)=\int_{A}q_{\eta}(x,\tl x)\dif\tl x,
\end{align}
for any $x\in\R^d$ and $A\in\mathcal B(\R^d)$.

The RLMC has  been shown to be an optimal discretization in a suitable sence in Cao et al. \cite{ca}, and it has been developed further in \cite{al,al2,hy,yu}. In particular, He et al. \cite[Theorem 1]{hy} showed that  \eqref{rmm} admits a unique stationary distribution $\pi_\eta$ and for any initial state $X_0\in\R^d$, the law of RLMC converges to $\pi_\eta$ under total variation distance,
but the authors did not give the convergence rate. Moreover, He et al. \cite[Proposition 2.2]{hy} derived an upper bound of Wasserstein-2 distance between $\pi$ and $\pi_\eta$ with order $O(\sqrt{\eta})$, which means that the RLMC would not converge to the real target $\pi$ and there is a  systematic error.

In this paper, we first establish the exponential ergodicity of RLMC (see Theorem \ref{ms} below) by verifying the 
irreducibility, strong Feller property and Lyapunov condition. Then, we design an decreasing step-size RLMC to approximate the target distribution $\pi$  and provide the convergence rate in terms of a functional class distance (see Theorem \ref{ms2} below).
 
 The rest of the paper is organized as follows. In Section \ref{mr}, we state our main results. The exponential ergodicity property of RLMC is proved in Section \ref{1021}, and Theorem \ref{ms2} is proved in Section \ref{pfms2}.

\vspace{10pt}
\textbf{Notations} We end this section by introducing some notations, which will be frequently used in the sequel. The inner product of $x, y \in\R^d$ is denoted by $\langle x,y\rangle$ and the Euclidean metric is denoted by $|x|$.
For any $u\in\mathbb{R}^{d}$ and matrix $A=(A_{ij})_{d\times d}$,  denote $Au^{\otimes 2}=\langle Au, u\rangle$, and the operator norm of the matrix $A$ is denoted by
$\|A\|_{\rm op}=\sup_{|u|=1}|Au|$. For any $k\geq 1$, let $\mathcal{C}^k(\R^d,\R)$ denote the collection of all  $k$-th order continuously differentiable functions from $\R^d$ to $\R$.  For any $f\in \mathcal{C}^k(\R^d,\R)$, let  $\nabla f$ and $\nabla^2f$ denote the gradient and the Hessian matrix of $f$, respectively. We define the supremum norm of $\|\nabla^if\|$ as
\begin{align*}
\|\nabla^if\|=\sup_{x\in\R^d}\|\nabla^if(x)\|_{\rm op}, \quad  i=1,2.
\end{align*}

Let $(\Omega, \mathcal F, \PP)$ be the probability space, for any random vector $X$, the law of $X$ is denoted by $\mathcal L(X)$ and the $L^2$-norm is defined by $\|X\|_2=[\E(|X|^2)]^{\frac{1}{2}}$. Let $\mathcal P(\R^d)$ be the space of probability measures on $\R^d$. Denote by $Leb(\cdot)$ the Lebesgue measure on $\R^d$.

 Let $V:\R^d\to \R_+$ be a measurable function. We define the  $V$-weighted supremum norm (\cite[p. 57]{haa}) of a function $h:\R^d\to \R$ as
\begin{align}
    \|h\|_V=\sup_{x\in\R^d}\frac{|h(x)|}{1+V(x)}.
\end{align}
The weighted total variation distance between $\mu,\nu\in\mathcal P(\R^d)$ is defined as
\begin{align}\label{wtv}
\dif_{\rm TV, V}(\mu,\nu)=\sup_{h: \|h\|_V\leq1}\left\{\int_{\R^d}h(x)\mu(\dif x)-\int_{\R^d}h(x)\nu(\dif x)\right\}.
\end{align}

Consider the following  function class
\begin{align*}
\mathcal G=\{h\in\mathcal{C}^2(\R^d,\R):\|\nabla^ih\|<+\infty \ {\rm with}\  i=1,2\}.
\end{align*}
Then for any $h\in\Gcal$, there exists a  positive constant $C_h$ such that $|h(x)|\leq C_h(1+|x|^2)$. Denote two probability measure spaces
\begin{align*}
\mathcal P_0=\left\{\mu\in\mathcal P(\R^d):\int_{\R^d}h(x)\mu(\dif x)<+\infty \  {\rm for\ any }\ h\in\Gcal\right\}
\end{align*}
and 
\begin{align*}
\mathcal P_2=\left\{\mu\in\mathcal P(\R^d):\int_{\R^d}|x|^2\mu(\dif x)<+\infty\right\},
\end{align*}
then it is easy to see that $\mathcal P_2\subset\mathcal P_0$. For $\mu_1,\mu_2\in\mathcal P_0$, define
\begin{align*}
\dif_{\Gcal}(\mu_1,\mu_2)=\sup_{h\in\Gcal}\left\{\int_{\R^d}h(x)\mu_1(\dif x)-\int_{\R^d}h(x)\mu_2(\dif x)\right\}.
\end{align*}
Obviously, $\dif_{\Gcal}(\cdot,\cdot)$ is a metric on $\mathcal P_0$.

\section{main results}\label{mr}
Let us begin with the following assumption on the potential $U$.
\begin{assumption}\label{ass}
(i) The potential function $U\in\mathcal{C}^2(\R^d,\R) $ and there exist two positive constant $0<m<L$ such that
\begin{align}\label{m}
m{\rm I}_{\rm d}\preceq \nabla^2U\preceq L{\rm I}_{\rm d}.
\end{align}

(ii) 
$\nabla U(0)=0$.

\end{assumption}
\begin{remark}
It follows from \cite[Lemma 4, Lemma 5]{dw} that \eqref{m} means $U$ is $m$-strongly convex and $L$-gradient Lipschitz  and it  is equivalent to
\begin{align*}
|\nabla U(x)-\nabla U(y)|\leq L|x-y|
\end{align*}
and 
\begin{align*}
\langle -\nabla U(x)+\nabla U(y),x-y\rangle \leq -m|x-y|^2.
\end{align*}
Taking $y=0$ implies that for any $x\in\R^d$,
\begin{align}\label{u1}
|\nabla U(x)|^2\leq L^2|x|^2
\end{align}
and
\begin{align}\label{u2}
\langle -\nabla U(x),x\rangle \leq -m|x|^2.
\end{align}
When $U$ is strongly convex, \cite[Theorem 2.1.8]{ne} shows that $U$ has a unique minimizer $x^*\in\R^d$ satisfying $\nabla U(x^*)=0$. In  Assumption \ref{ass} (ii),   we assume $x^*=0$, which is merely for computational convenience, the fact that  $\nabla U(0)$ is not zero does not affect the main results of this paper.
\end{remark}
Our first result is the following exponentially ergodicity of RLMC \eqref{rmm}.
\begin{theorem}\label{ms}
Under Assumption \ref{ass}, for any  step size $\eta\in(0,\frac{m}{L^2}]$, the Markov chain $(X_{k\eta})_{k\in\mathbb N_0}$ is exponentially ergodic with unique invariant measure $\pi_\eta$.  More precisely, let
$V(x)=1+|x|^2$ and $\nu Q_\eta^n$ be the law of $X_{n\eta}$ with initial distribution $\nu$ satisfying $\nu(V)<\infty$, then
\begin{eqnarray*}
d_{\rm TV,V} (\nu Q_\eta^n,\pi_\eta)&\leq& C(1+\nu(V))e^{-cn},
\end{eqnarray*}
where $C$ and $c$ are some positive constants independent of $n$.
\end{theorem}

\begin{remark}
     He et al. \cite[Theorem 1]{hy} showed that, for every $x\in\R^d$, $\delta_x Q_\eta^n$ converges to $\pi_\eta$ under total variation distance, but the authors did not give any information about the convergence rate. Note that for $\mu,\nu\in\mathcal P(\R^d)$, the total variation distance is defined as
     \begin{align*}
         \dif_{\rm TV}(\mu, \nu)=\sup_{h: \|h\|_{\infty}\leq1}\left\{\int_{\R^d}h(x)\mu(\dif x)-\int_{\R^d}h(x)\nu(\dif x)\right\}.
     \end{align*}
     Compared with \eqref{wtv}, one has the relation $\dif_{\rm TV}(\mu, \nu)\leq\dif_{\rm TV,V}(\mu, \nu)$. In Theorem \ref{ms}, we have established the convergence rate under the distance $d_{\rm TV,V}$,  which  is an improvement of  He et al. \cite[Theorem 1]{hy}.
\end{remark}

It has been proved in \cite{hy} that choosing a constant step-size for RLMC leads to sampling bias.
We consider the  random midpoint  Euler-Maruyama scheme  with decreasing step-size sequence $\Gamma=(\gamma_n)_{n\in \mathbb{N}}$ below.
Let  $\gamma_k>0$ the $k$-th step size, $t_0=0$, $t_k=\sum_{i=1}^k \gamma_i$, $k \geq 1$. Denote
\begin{align}\label{bbc}
Y_{t_{k+1}}=Y_{t_{k}}-\gamma_{k+1}\nabla U(Y_{t_{k}+u_{k+1}\gamma_{k+1}})+\sqrt{2\gamma_{k+1}}\xi_{k+1},
\end{align}
where the random midpoint $Y_{t_{k}+u_{k+1}\gamma_{k+1}}$ is given by
\begin{align*}
Y_{t_{k}+u_{k+1}\gamma_{k+1}}=Y_{t_{k}}-u_{k+1}\gamma_{k+1}\nabla U(Y_{t_{k}})+\sqrt{2u_{k+1}\gamma_{k+1}}\xi_{k+1}^\prime.
\end{align*}
Here, the setting of $\{u_k\}_{k\geq1}$, $(\xi_k)_{k\geq1}$ and $(\xi^\prime_k)_{k\geq1}$ are the same as that of \eqref{rmm}.

We make the following assumption on the step sequence $\Gamma=(\gamma_n)_{n\in\mathbb{N}}$.
\begin{assumption}\label{assump2}
(i) Let $(\gamma_n)_{n\in \mathbb{N}}$ be a sequence of positive and non-increasing step sizes satisfying $\lim\limits_{n\to\infty}\gamma_n=0$ and  $\sum\limits_{n\geq1}\gamma_n=+\infty;$

(ii)
\begin{align*}
\omega:=\limsup_{n\rightarrow\infty}\frac{\gamma_{n}^{1/2}-\gamma_{n+1}^{1/2}}{\gamma_{n+1}^{3/2}}<+\infty.
\end{align*}
\end{assumption}

We next  show that the marginal distribution of the $Y_{t_n}$   converges to the target distribution $\pi$ and provide explicit convergence bounds under the $d_{\mathcal{G}}$ metric. 

\begin{theorem}\label{ms2}

Let the potential  $U$  satisfy Assumption \ref{ass}, and  $(X_t^x)_{t\geq0}$ be the Langevin diffusion with unique  stationary distribution $\pi$ and initial value $x$. Assume that the step-size sequence $\Gamma$ satisfy Assumption \ref{assump2} with $\omega<m/2$, then we have
\begin{align}\label{convergence1}
d_{\mathcal{G}} \left(\mathcal L\big( X_{t_n}^{x}\big),\mathcal L\big( Y_{t_n}^{x} \big)\right) \leq C (1+|x|^2)\gamma_n^{1/2}.
\end{align}
Furthermore, there exists some positive constant $C$ independent of $\Gamma$ such that
\begin{align}
d_{\mathcal{G}} \left(\mathcal L\big(Y_{t_n}^{x} \big), \pi\right) \leq C (1+|x|^2) \gamma_n^{1/2}.
\end{align}\label{convergence2}
\end{theorem}

\section{Exponential ergodicity  of RLMC}\label{1021}
In this section, we prove Theorem \ref{ms} by verifying that the Markov process $(X_{k\eta})_{k\in\mathbb N_0}$ defined by \eqref{rmm} is irreducible, strongly Feller, and satisfies the Lyapunov condition. These properties will be stated as lemmas in the Subsection \ref{e}, and then we provide a detailed proof of the main theorem based on these lemmas. The  proofs of  lemmas will be presented in the Subsection \ref{f}.

\subsection{Proof of Theorem \ref{ms}}\label{e}
We recall the definition of the irreducibility and strong Feller property of a Markov chain, which could be found in \cite{MC} and \cite{ha}.

\begin{definition}[Accessible set, small set, irreducibility]\label{small}
Let $P$ be a Markov transition kernel on the state space $(X,\mathcal X)$.

(i) A set $A\in \mathcal X$ is said to be accessible if for all $x\in X$, there exists an integer $n\geq1$ such that $P^n(x,A)>0$.

(ii) A set $A\in \mathcal X$ is called a small set if there exists a positive integer $n$ and a nonzero measure $\mu$ on $(X,\mathcal X)$ such that for all $x\in A$ and $B\in\mathcal X$,
\begin{align*}
P^n(x,B)\geq\mu(B).
\end{align*}
The set $A$ is then said to be an $(n,\mu)$-small set.

(iii) $P$ is said to be irreducible if it admits an accessible small set.
\end{definition}

\begin{definition}[Hairer \cite{ha}]\label{sf}
A Markov kernel $P$ on the state space $(X,\mathcal X)$ has the strong Feller property if $Pf$ is continuous for every bounded measurable function $f:X\to \mathbb{R}$.
\end{definition}

It turns out that the transition kernel $Q_\eta$ given by \eqref{zyh} is irreducible and strong Feller.

\begin{lemma}\label{irr}
Under Assumption \ref{ass}, for any step size $\eta\in(0,1)$,  the transition  kernel $Q_\eta$  of the
Markov chain $(X_{k\eta})_{k\in\mathbb N_0}$ defined by \eqref{rmm} is irreducible.
\end{lemma}

\begin{lemma}\label{sfe}
  Under Assumption \ref{ass}, for any  step size $\eta\in(0,1]$,  $Q_\eta$ is strong Feller.
\end{lemma}

Furthermore, the following Lyapunov function condition holds  for $Q_\eta$.
\begin{lemma}\label{lya}
Let $V(x)=|x|^2$, $x\in\R^d$. Under Assumption \ref{ass}, for any  step size $\eta\in(0,\frac{m}{L^2}]$,   the Markov kernel $Q_\eta$ satisfies
\begin{align*}
Q_\eta V(x)\leq \lambda V(x)+b,
\end{align*}
where $\lambda\in(0,1)$ and  $b>0$ are two constants depending on $\eta, m, L$ and $d$.
\end{lemma}

With the above auxiliary lemmas in hand, we may get the exponential ergodicity of $(X_{k\eta})_{k\in\mathbb N_0}$.

\begin{proof}[Proof of Theorem \ref{ms}]
According to Meyn and Tweedie \cite[Theorem 6.1]{me}, $(X_{k\eta})_{k\in\mathbb N_0}$ is exponentially ergodic with ergodic measure $\pi_\eta$, and  we can get the following inequality:
\begin{eqnarray}\label{general}
\dif_{\rm TV,V}(\delta_{x}Q_\eta^n,\pi_\eta)\leq C(1+V(x)) e^{-cn},
\end{eqnarray}
where $\delta_{x}$ is the delta measure on $x$. Then, for any probability measure $\nu$ with $\nu(V) <\infty$, \eqref{general} implies that
\begin{align*}
\dif_{\rm TV,V}(\nu Q_\eta^n,\pi_\eta)=&
\ \dif_{\rm TV,V}(\int_{\R^d}\delta_x Q_\eta^n\nu(\dif x),\pi_\eta)\\
=&\ \sup_{f:\|f\|_V\leq1}\left\{\int_{\R^d}[\delta_x Q_\eta^n(f)-\pi_\eta(f)]\nu(\dif x)\right\}\\
=&\ \int_{\R^d}\dif_{\rm TV,V}(\delta_{x}Q_\eta^n,\pi_\eta)\nu(\dif x)\\
\leq &\ C(1+\nu(V)) e^{-cn}.
\end{align*}

\end{proof}

\subsection{Proof of Lemmas   }\label{f}

\begin{proof}[Proof of Lemma \ref{irr}]
We claim that for any step size $\eta\in(0,1)$,
 if $E$ is a compact subset of $\R^d$ such that  $Leb(E)>0$, then $E$ is an accessible $(1,\epsilon\nu)$- small set, where
\begin{align*}
\epsilon=\inf_{x\in C,\tl x\in C}q_\eta(x,\tl x),
\end{align*}
$q_\eta$ is the transition density function defined in \eqref{zymd} and $\nu(\cdot)=Leb(\cdot\cap E)$. Indeed,
suppose that $E$ is a compact subset of $\R^d$ such that $Leb(E)>0$, then for all $x\in E$ and $A\in\mathcal B(\R^d)$,
\begin{align}\label{436}
Q_\eta(x,A)=\int_Aq_\eta(x,\tl x)\dif \tl x\geq \int_{A\cap E}q_\eta(x,\tl x)\dif y\geq \inf_{x\in E,\tl x\in E}q_\eta(x,\tl x) Leb(A\cap E).
\end{align}
It is easy to check that the transition density function $q_{\eta}(x,\tl x)$ in \eqref{zymd} is a bivariate continuous function. So there exist a point $(x_0,y_0)$ belonging to  $E$ such that
 \begin{align}\label{437}
\inf_{x\in E,\tl x\in E}q_\eta(x,\tl x)=q_\eta(x_0,y_0)\in(0,1].
\end{align}
Combining \eqref{436} and \eqref{437}, we have for all $x\in E$ and $A\in\mathcal B(\R^d)$
\begin{align*}
Q_\eta(x,A)\geq\epsilon \nu(A),
\end{align*}
where $\nu(\cdot)=Leb(\cdot\cap E)$ and
\begin{align*}
\epsilon=\inf_{x\in E,\tl x\in E}q_\eta(x,\tl x)\in(0,1].
\end{align*}
Particularly, $Q_\eta(x,E)\geq\epsilon Leb (E)>0$ for all $x\in\R^d$.
By Definition \ref{small}, $E$ is a accessible $(1,\epsilon\nu)$-small set, and then $Q_\eta$ is irreducible.

\end{proof}

\begin{proof}[Proof of Lemma \ref{sfe}]
If $f$ is a bounded measurable function on $\R^d$, then clearly $Q_\eta f$ is bounded. It remains to verify that $Q_\eta f$ is a continuous funtion, i.e. for any $x\in\R^d$, if $x_n\to x$, then $Q_\eta f(x_n)\to Q_\eta f(x)$.

  If $x_n\to x$, there exists a positive constant $M$ such that $|x_n|\leq M$  for all $n\geq1$.  Recall the transition density $q_\eta$ in \eqref{zymd}, we have
 \begin{align*}
 Q_\eta f(x_n)=\int_{\R^d}f(y)q_\eta(x_n,y)\dif y.
 \end{align*}
  Since $q_\eta$ is  continuous and $x_n\to x$, for any $y\in \R^d$,
  \begin{align*}
 f(y)q_\eta(x_n,y)\to f(y)q_\eta(x,y),\quad n\to\infty.
 \end{align*}
  Notice that
  \begin{align*}
 f(y)q_\eta(x_n,y)\leq \sup_{x:|x|\leq M}f(y)q_\eta(x,y).
 \end{align*}
Since the right-hand side is integrable, by the dominated convergence theorem,
\begin{align*}
 Q_\eta f(x_n)=\int_{\R^d}f(y)q_\eta(x_n,y)\dif y\to\int_{\R^d}f(y)q_\eta(x,y)\dif y=Q_\eta f(x), \quad n\to\infty.
 \end{align*}
 The proof is complete.
\end{proof}

\begin{proof}[Proof of Lemma \ref{lya}]
Recall the definition of the Markov chain $(X_{k\eta})_{k\geq0}$ in \eqref{rmm}, for any initial state $X_0=x\in\R^d$, let
\begin{align}\label{one-step}
X_{u\eta}=x-u\eta\nabla U(x)+\sqrt{2u\eta}\xi^\prime\\
X_1=x-\eta\nabla U(X_{u\eta})+\sqrt{2\eta}\xi,
\end{align}
where $u\thicksim U[0,1]$ and $\xi,\xi^\prime$ are two independent standard normal random vectors.

Then 
\begin{align}\label{1955}
Q_\eta V(x)&=\E[|X_1|^2]\nonumber\\
&=\E[|x-\eta\nabla U(X_{u\eta})+\sqrt{2\eta}\xi|^2]\nonumber\\
&=\E\Big[|x|^2+\eta^2|\nabla U(X_{u\eta})|^2+2\eta d+2\langle x,-\eta\nabla U(X_{u\eta})\rangle\Big],
\end{align}
where the last equation holds since $\xi$ is a centered $d$-dimensional normal random vector independent of $X_{u\eta}$.

 Next we focus on  the terms $\E[|\nabla U(X_{u\eta})|^2]$ and $\E[2\langle x,-\eta\nabla U(X_{u\eta})\rangle]$.
First, by \eqref{u1}, \eqref{u2}, \eqref{one-step} and the independence of $u$ and $\xi^\prime$, we have
\begin{align}\label{1438}
\E[|\nabla U(X_{u\eta})|^2]\leq&L^2\E[|X_{u\eta}|^2]\nonumber\\
=&L^2\E\Big[|x|^2+u^2\eta^2|\nabla U(x)|^2+2u\eta |\xi^\prime|^2+2\langle x,-u\eta\nabla U(x)\rangle\Big]\nonumber\\
\leq&L^2\Big(|x|^2+\frac{1}{3}\eta^2L^2|x|^2+\eta d-m\eta|x|^2\Big)\nonumber\\
\leq&\left(L^2+\eta^2L^4-m\eta L^2\right)|x|^2+L^2\eta d.
\end
{align}
By Taylor's formula, there exists a random point $x_0$ depending on $x$ and $X_{u\eta}$ such that 
\begin{align*}
\nabla U(X_{u\eta})=\nabla U(x)+\nabla^2U(x_0)(X_{u\eta}-x).
\end{align*}
Hence by \eqref{m}, \eqref{u1}, \eqref{u2} and the independence of $u$ and $\xi^\prime$, we have
\begin{align}\label{1439}
&\E[2\langle x,-\eta\nabla U(X_{u\eta})]\nonumber\\
=&\E\Big[\langle 2\langle x,-\eta\nabla U(x)\rangle +2\langle \eta x,u\eta\nabla^2U(x_0)\nabla U(x)\rangle+2\langle \eta x,-\sqrt{2u}\eta^{\frac{1}{2}}\nabla^2U(x_0)\xi^\prime\rangle\Big]\nonumber\\
\leq&-2m\eta|x|^2+L^2\eta^2|x|^2+\eta^2|x|^2+\eta L^2d
\end{align}
Combining \eqref{1955}-\eqref{1439} and assuming $\eta\leq \frac{m}{L^2}$ , we obtain
\begin{align}\label{1230}
Q_\eta V(x)&\leq\left[1-2m\eta+(2L^2+1)\eta^2\right]|x|^2+(2+L^2)d\eta +L^2d\eta^3.
\end{align}
The proof is complete.
\end{proof}

\section{Convergence of decreasing step-size RLMC }\label{pfms2}
In this section, we will  prove Theorem \ref{ms2}. The proof mostly rely on the method developted in Pag\`es and Panloup \cite{pages}, which has been widely used in recent years to study the Euler-Maruyama approximations of ergodic stochastic differential equations, such as \cite{ch,  jtp2, de,  gu}. 
In Subsection \ref{sec:result}, we provide the necessary preliminaries and a series of lemmas that will be used to establish the main results. We prove the main theorem by combining these lemmas with the approach  in \cite{pages} in Subsection \ref{Last}. To keep the structure concise, all proofs of the lemmas are deferred to Subsection \ref{1023}.

\subsection{Preliminaries}\label{sec:result}

Denote by $X_t^{x}$ the solution to SDE in \eqref{ld}, where the superscript $x$ indicate the dependence upon the initial value $X_0=x$. For any $0\leq s\leq t<\infty$, denote the semigroup of $(X_t)_{t\geq 0}$ by $P_{s,t}$, that is, for any bounded Borel function $f\in \mathcal{B}_b(\R^d,\R)$ and any $x\in \R^d $, one has
\begin{eqnarray*}
P_{s,t}f(x) = \E\left[f(X_t)\Big|X_{s}=x\right].
\end{eqnarray*}
Notice that the process $(X_t)_{t\geq 0}$ is time homogeneous, hence we  denote
\begin{align*}
P_{t-s}f(x)=P_{s,t}f(x).
\end{align*}
The discrete time 
 process
$(Y_{t_k})_{k\in \mathbb{N}_0}$ given in \eqref{bbc} is inhomogeneous, we denote
\begin{align*}
\tilde{P}_{t_{i},t_{j}}f(x)=\mathbb{E}\left[f(Y_{t_{j}})\Big|Y_{t_{i}}=x\right], \quad 0\leq i\leq j<\infty, f\in\mathcal{B}_b(\R^d,\R)
\end{align*}
and the probability trasition
\begin{align*}
\tilde{P}_{\gamma_{k+1}} f(x)=\tilde{P}_{t_{k},t_{k+1}}f(x), \quad k\in\mathbb{N}_0.
\end{align*}
Let $(\tilde{Y}_{\gamma}^{x})$ denote the one step random midpoint Euler-Maruyama scheme with any step size $\gamma\in(0,\gamma_{1}]$ starting from $x$ at time 0, that is,
\begin{align}\label{onestep}
\tl Y_{\gamma}^x=x-\gamma\nabla U(Y_{u\gamma})+\sqrt{2\gamma}\xi,
\end{align}
where 
\begin{align*}
Y_{u\gamma}=x-u\gamma\nabla U(x)+\sqrt{2u\gamma}\xi^\prime,
\end{align*}
and then
\begin{eqnarray*}
\E f(\tilde{Y}_{\gamma}^{x})=\tilde{P}_{\gamma} f(x).
\end{eqnarray*}

First, we  give the moment estimates of   $(X_t^x)_{t\geq0}$ and  $(Y_{t_{k}}^x)_{k\geq 0}$ .

\begin{lemma}\label{LYexp}
Let $\left(X_{t}^{x}\right)_{t\geq0}$ be the Langevin diffusion in \eqref{ld} and $(Y^x_{t_{k}})_{k\geq 0}$ be the RLMC with decreasing step size  defined in \eqref{bbc}. Then under Assumptions \ref{ass},  for any initial state $x$ and $t\geq 0$, there exists a constant $C>0$ independent of $t$ such that
\begin{eqnarray}\label{Xmoment}
\mathbb{E}|X_{t}^{x}|^{2}\leq C(1+|x|^2),
\end{eqnarray}
and for any $k\geq1$, there exists $C>0$ independent of $\Gamma=(\gamma_n)_{n\in \mathbb{N}}$ such that
\begin{eqnarray}\label{Yexp}
\mathbb{E}|Y_{t_{k}}^{x}|^{2}\leq C(1+|x|^{2}).
\end{eqnarray}
\end{lemma}

Next we state the  following  propositions which are crucial in the proof of our main result.
\begin{lemma}\label{tidu}
Let $(P_t)_{t\geq0}$ be the semigroup of Langevin diffusion $(X_t)_{t\geq0}$. Under Assumption \ref{ass}, there exists some positive constant $C$ independent of $t$ such that for any $h\in\Gcal$ and $t\geq0$
\begin{align*}
\|\nabla P_th\| &\leq \| \nabla h \|e^{-mt},  \\
\|\nabla^2 P_th\|&\leq C (\| \nabla h \|+\| \nabla^2 h \|) e^{-\frac{m}{2}t}
\end{align*}
 where $m$ is the constant in \eqref{m}.
\end{lemma}

\begin{lemma}\label{prop:onestep}
 Under Assumptions \ref{ass}, for any $x\in\R^d$, $0<\gamma\leq \gamma_1<1$ and $f\in\mathcal{C}^2(\R^d,\R)$,  there exists some positive constant $C$  independent of $\gamma$ such that 
\begin{align*}
|\E f(X_{\gamma}^{x})-\E f(\tilde{Y}_{\gamma}^{x})|\leq C(1+|x|^2)(\|\nabla f\|+\|\nabla^{2} f\|)\gamma^{\frac{3}{2}}.
\end{align*}
\end{lemma}

\begin{lemma}\label{Step}
Let $(\gamma_n)_{n\in \mathbb{N}}$ be a non-increasing positive sequence. If
\begin{align*}
\omega=\limsup_{n\rightarrow\infty} \{ (\gamma_{n}^{1/2}-\gamma_{n+1}^{1/2}) \gamma_{n+1}^{-3/2} \}<+\infty,
\end{align*}
let $m>2\omega$ and $(u_n)_{n\in \mathbb{N}}$ be a sequence defined by $u_0=0$ and
\begin{align*}
u_n=\sum_{k=1}^n \gamma_k^{3/2}e^{-\frac{m}{2}(t_n-t_k)}, \ \ \   n\geq1.
\end{align*}
Then, we have
\begin{align}\label{s1}
\limsup_{n\rightarrow\infty} \{ u_n \gamma_n^{-1/2} \}<+\infty,
\end{align}
and 
\begin{align}\label{exp1}
\lim_{n\rightarrow\infty} e^{-m t_{n}} \gamma_{n}^{-\frac{1}{2}}=0.
\end{align}
\end{lemma}

\subsection{Proof of Theorem \ref{ms2}}\label{Last}

With the help of the above preliminary propositions and lemmas, we are in the position to prove Theorem \ref{ms2}.

\begin{proof}[Proof of Theorem \ref{ms2}]
For any $n\geq1$, by the definition of the  $d_\Gcal$, we have
\begin{align}\label{R1}
d_\Gcal \left(\Lcal\left( X_{t_n}^{x}\right),\Lcal\left( Y_{t_n}^{x}\right)\right) 
=&\sup_{h\in\Gcal}\{\E h(X_{t_n}^{x})-\E h(Y_{t_n}^{x}) \}\nonumber\\
=&\sup_{h\in\Gcal} \{P_{\gamma_1}\circ \cdots\circ P_{\gamma_n}h(x)-\tl{P}_{\gamma_1}\circ \cdots\circ\tl{P}_{\gamma_n}h(x) \} \nonumber\\
\leq&\sup_{h\in\Gcal}\left\{ \sum_{k=1}^n|\tl{P}_{\gamma_1}\circ \cdots\circ\tl{P}_{\gamma_{k-1}}\circ(P_{\gamma_{k}}-\tl{P}_{\gamma_k})\circ P_{t_n-t_k}h(x)| \right\}.
\end{align}
It follows from Lemma \ref{tidu} and \ref{prop:onestep} that for any $k \in\{1,2,\cdots,n\}$,
\begin{align*}
|(P_{\gamma_{k}}-\tl{P}_{\gamma_k})\circ P_{t_n-t_k}h(x)|
\leq& C(1+|x|^2)(\|\nabla P_{t_n-t_k}h\|+\|\nabla^2 P_{t_n-t_k}h\|)\gamma_k^{\frac{3}{2}}\\
\leq& C(1+|x|^2)(\|\nabla h\|+\|\nabla^2h\|)e^{-\frac{m}{2}(t_n-t_k)}\gamma_k^{\frac{3}{2}}.
\end{align*}
Then, integrating with respect to $\tl{P}_{\gamma_1}\circ \cdots\circ\tl{P}_{\gamma_{k-1}}$ yields
\begin{align}\label{R2}
&|\tl{P}_{\gamma_1}\circ \cdots\circ\tl{P}_{\gamma_{k-1}}\circ(P_{\gamma_{k}}-\tl{P}_{\gamma_k})\circ P_{t_n-t_k}h(x)|\nonumber \\
\leq&  C(\|\nabla h\|+\|\nabla^2h\|)e^{-\frac{m}{2}(t_n-t_k)}\gamma_k^{\frac{3}{2}}(1+\E|Y_{t_{k-1}}^{x}|^2)\nonumber\\
\leq&  C(1+|x|^2)(\|\nabla h\|+\|\nabla^2h\|)e^{-\frac{m}{2}(t_n-t_k)}\gamma_k^{\frac{3}{2}},
\end{align}
where the last inequality holds by (\ref{Yexp}). Hence, combining \eqref{s1}, \eqref{R1} and \eqref{R2} yields
\begin{align}\label{dh0}
d_\Gcal \left(\Lcal\left( X_{t_n}^{x}\right),\Lcal\left( Y_{t_n}^{x}\right)\right) &\leq  C(1+|x|^2)\sum_{k=1}^ne^{-\frac{m}{2}(t_n-t_k)}\gamma_k^{\frac{3}{2}}
\leq C_1(1+|x|^2)\gamma_n^{\frac{1}{2}}.
\end{align}

Next we consider $d_\Gcal \left(\Lcal\left( X_{t_n}^{x}\right),\pi\right)$.  According to 
Durmus and Moulines \cite[Proposition 1 ]{du2}, under Assumption \ref{ass},
\begin{align*}
W_2\left(\Lcal\left( X_{t_n}^{x}\right),\pi\right)\leq e^{-mt_n}\{|x|+(d/m)^{1/2}\}.
\end{align*}

Denote
\begin{align*}
\mathcal F =\{f|f:\R^d\to\R,\|f\|_{\rm Lip}\leq 1\}.
\end{align*}
By the definition of the function class $\mathcal G$,  there exists a constant $M_1$   such that $\|\nabla h\|\leq M_1$ for any $h\in\mathcal G$. So we have
\begin{align}\label{dh1}
d_\Gcal \left(\Lcal\left( X_{t_n}^{x}\right),\pi\right)=\sup_{h\in\mathcal G}M_1\{\E \frac{h}{M_1}(X_{t_n}^{x})-\pi(\frac{h}{M_1})\}\leq M_1\sup_{h\in\mathcal F}\{\E h(X_{t_n}^{x})-\pi(h)\}.
\end{align}

 By Kantorovich dual formula \cite[Remark 6.5]{v}, we have
 \begin{align}\label{dh2}
&\sup_{h\in\mathcal F}\{\E h(X_{t_n}^{x})-\pi(h)\}\nonumber\\
=&W_1\left(\Lcal\left( X_{t_n}^{x}\right),\pi\right)\leq W_2\left(\Lcal\left( X_{t_n}^{x}\right),\pi\right)\leq e^{-mt_n}\{|x|+(d/m)^{1/2}\}.
\end{align}

Furthermore, combining \eqref{exp1}, \eqref{dh0}, \eqref{dh1} and  \eqref{dh2}  together yields 
\begin{align*}
 d_\Gcal \left(\Lcal\left( Y_{t_n}^{x}\right),\pi\right)
 \leq& d_\Gcal \left(\Lcal\left( X_{t_n}^{x}\right),\Lcal\left( Y_{t_n}^{x}\right)\right)+d_\Gcal \left(\Lcal\left( X_{t_n}^{x}\right),\pi\right)\\
 \leq&C_1(1+|x|^2)\gamma_n^{\frac{1}{2}}+M_1\{|x|+(d/m)^{1/2}\}e^{-mt_n}\\
 \leq& C\left(1+|x|^{2}\right)\gamma_{n}^{\frac{1}{2}}.
\end{align*}

The proof is complete.
\end{proof}

\subsection{Proofs of Lemmas}\label{1023}
\begin{proof}[Proof of Lemma \ref{LYexp}]
(i) By  It\^{o}'s formula,
\begin{align*}
|X_t^x|^2=|x|^2+\int_0^t\langle 2X_s^x, -\nabla U(X_s^x)\rangle\dif s+\int_0^t\langle 2X_s^x, \sqrt2 \dif B_s\rangle+2dt.
\end{align*}
Therefore by \eqref{u2} we have
\begin{align*}
\frac{\dif}{\dif t}\E|X_t^x|^2=&\E\langle 2X_t^x, -\nabla U(X_t^x)\rangle+2d\\
\leq&-2m\E|X_t^x|^2+2d,
\end{align*}
which implies that
\begin{align*}
\E|X_t^x|^2\leq |x|^2e^{-2mt}+\frac{d}{m}(1-e^{-2mt}).
\end{align*}

(ii) By \eqref{1230}, there exist positive constants $C_{1},C_{2}$ such that
\begin{align*}
\mathbb{E}|Y_{t_{k+1}}|^{2}\leq&\left(1-2m\gamma_{k+1}+C_{1}\gamma_{k+1}^{2}\right)\mathbb{E}|Y_{t_{k}}|^{2}+C_{2}\gamma_{k+1}.
\end{align*}
Since $\lim_{k\rightarrow\infty}\gamma_{k}=0$, there exists $k_{0}\in\mathbb{N}$, such that, for every  $k\geq k_{0}$ satisfying $C_{1}\gamma_{k+1}^{2} \leq m\gamma_{k+1}$ and
 \begin{align*}
\mathbb{E}|Y_{t_{k+1}}|^{2}\leq&\left(1-m\gamma_{k+1}\right)\mathbb{E}|Y_{t_{k}}|^{2}+C_{2}\gamma_{k+1},
\end{align*}
and $1-m\gamma_{k+1}>0$. Inductively, for any $k>k_{0}$, we have
\begin{align*}
\mathbb{E}|Y_{t_{k}}|^{2}\leq\prod_{_{j=k_{0}+1}}^{k}\left(1-m\gamma_{j}\right)\mathbb{E}|Y_{t_{k_{0}}}|^{2}
+C_{2}\sum_{j=k_{0}+1}^{k}\left[\gamma_{j}\prod_{l=j+1}^{k}\left(1-m\gamma_{l}\right)\right],
\end{align*}
here, we have adopted the convention that $\prod_{l=k+1}^{k}\left(1-m\gamma_{l}\right)=1$. 
Note that
\begin{align*}
m\sum_{j=k_{0}+1}^{k}\gamma_{j}\prod_{l=j+1}^{k}\left(1-m\gamma_{l}\right)=1-\prod_{l=k_{0}+1}^{k}\left(1-m\gamma_{l}\right)\leq1.
\end{align*}
Hence, the above two inequalities imply
\begin{align}\label{induc1}
\mathbb{E}|Y_{t_{k}}|^{2}\leq\mathbb{E}|Y_{t_{k_{0}}}|^{2}
+\frac{C_{2}}{m}, \quad \forall k\geq k_0.
\end{align}
In addition, by a standard argument (see Lamberton and Pag\`es \cite[Lemma 2]{LP02}), it is easy to verify that for every $0\leq k\leq k_{0}$,
\begin{align}\label{induc2}
\mathbb{E}|Y_{t_{k}}^{x}|^{2}\leq C(1+|x|^{2}).
\end{align}
Combining (\ref{induc1}) and (\ref{induc2}), for any $k\geq0$, we have
\begin{align*}
\mathbb{E}|Y_{t_{k}}^{x}|^{2}\leq C(1+|x|^{2}).
\end{align*}
\end{proof}

\begin{proof}[Proof of Lemma \ref{tidu}]
The proof of this lemma is omitted. In fact, it can be proved by following the line of proof of Proposition 3.4 in \cite{jtp2} with minor modification, see also \cite{ss}.

\end{proof}

To prove Lemma \ref{prop:onestep}, we need the following additional lemma.
\begin{lemma}\label{one}
Let $\left(X_{t}^{x}\right)_{t\geq0}$ be the solution to SDE (\ref{ld}) and $(\tilde{Y}_{\gamma}^{x})$ be the one step random midpoint Euler-Maruyama  scheme given by \eqref{onestep} with step $\gamma\in(0,\gamma_{1}]$, $\gamma_1<1$. Under Assumptions \ref{ass},  for every $t\in[0,\gamma]$, there exists some positive constant $C$ independent of $t$ such that
\begin{align}
\|X_{t}^{x}-x\|_2&\leq C (1+|x|)t^{\frac{1}{2}},\label{Xdiff}\\
\|\tilde{Y}_{\gamma}^{x}-x\|_2&\leq C (1+|x|)\gamma^{\frac{1}{2}},\label{Ydiff}\\
\|X_{\gamma}^{x}-\tilde{Y}_{\gamma}^{x}\|_2&\leq C (1+|x|)\gamma^{\frac{3}{2}}\label{X-Ydiff}.
\end{align}
\end{lemma}

\begin{proof}
(i) Recall that
\begin{align*}
X_t^x=x-\int_0^t\nabla U(X_s)\dif s+\sqrt2B_t,
\end{align*}
it follows from Minkowski inequality, \eqref{u1} and  Lemma \ref{LYexp} that 
\begin{align*}
\|X_t^x-x\|_2=&\left\|\int_0^t\nabla U(X_s)\dif s\right\|_2+\sqrt{2t}\|B_1\|_2\\
&\leq L\int_0^t\left\|X_s\right\|_2\dif s+\sqrt{2t}\\
&\leq C(1+|x|)t^{\frac{1}{2}}.
\end{align*}

(ii) By the definition of $\tl Y_{\gamma}^x$ in \eqref{onestep}, for $0<\gamma<\gamma_1<1$
\begin{align}\label{eati-Yur}
\|\tl Y_{\gamma}^x-x\|_2\leq &\gamma\|\nabla U(Y_{u\gamma})\|_2+\sqrt{2\gamma}\|\xi\|_2\nonumber\\
\leq&\gamma L\|x-u\gamma\nabla U(x)+\sqrt{2u\gamma}\xi^\prime\|_2+\sqrt{2\gamma}\|\xi\|_2\nonumber\\
\leq& C(1+|x|)\gamma^{\frac{1}{2}},
\end{align}
where the last inequality is due to \eqref{u1} and the independence of $u$ and $\xi^\prime$.

(iii) By Minkowski inequality, 
\begin{align*}
\|X_{\gamma}^{x}-\tilde{Y}_{\gamma}^{x}\|_2&\leq
\gamma\|\nabla U(Y_{u\gamma})-\nabla U(x)\|_2+\Big\|\int_0^\gamma\nabla U(X_s)-\nabla U(x)\dif s\Big\|_{2}.
\end{align*}
By \eqref{u1} and the similar analysis as that in \eqref{eati-Yur}, it is easy to see that the first term is controlled by $C(1+|x|)\gamma^{\frac{3}{2}}. $
For the second term,
\begin{align*}
\Big\|\int_0^\gamma\nabla U(X_s)-\nabla U(x)\dif s\Big\|_{2}=&\sqrt{\E\left[\left|\int_0^\gamma\nabla U(X_s)-\nabla U(x)\dif s\right|^2\right]}\\
\leq&\gamma^{\frac{1}{2}}\sqrt{\int_0^\gamma\E[|\nabla U(X_s)-\nabla U(x)|^2]\dif s}\\
\leq&\gamma L\left\{\E\left[\sup_{0\leq t\leq \gamma}\|X_t-X_0\|^2\right]\right\}^{\frac{1}{2}}.
\end{align*}
According to  Lemma 2 in \cite{hy},
\begin{align*}
\E\left[\sup_{0\leq t\leq \gamma}\|X_t-X_0\|^2\right]\leq
4\gamma^2|\nabla U(x)|^2+8L^2d\gamma^3+2d\gamma.
\end{align*}
Combining the above inequalities and taking into account that $0<\gamma<1$, we obtain
\begin{align}\label{tsg2}
\Big\|\int_0^\gamma\nabla U(X_s)-\nabla U(x)\dif s\Big\|_{2}\leq C(1+|x|)\gamma^{\frac{3}{2}}.
\end{align}
The proof is complete.

\end{proof}

\begin{proof}[Proof of Lemma \ref{prop:onestep}]
 By  Taylor formula, for any $f\in \mathcal{C}^2(\R^d,\R)$ , $x,y,z\in\R^d$,
\begin{align*}
f(z)-f(y)=&\langle\nabla f(x),z-y\rangle+\int_0^1\langle \nabla^2f(x+u(y-x))(y-x),z-y\rangle\dif u\\
&+\int_0^1(1-u)\nabla^2 f(uz+(1-u)y)(z-y)^{\otimes 2}\dif u.
\end{align*}
Applying this expansion with $z=\tl{Y}_{\gamma}^{x}$ and $y=X_{\gamma}^{x}$  yields
\begin{align}\label{I0}
|\E f(\tl{Y}_{\gamma}^{x})-\E f(X_{\gamma}^{x})|=\mathcal{I}_1+\mathcal{I}_2+\mathcal{I}_3,
\end{align}
where
\begin{align*}
&\mathcal{I}_1=\E[|\langle\nabla f(x),\tl{Y}_{\gamma}^{x}-X_{\gamma}^{x}\rangle|],\\
&\mathcal{I}_2=\E\left[\left|\int_0^1\langle \nabla^2f(x+u(X_{\gamma}^{x}-x))(X_{\gamma}^{x}-x),\tl{Y}_{\gamma}^{x}-X_{\gamma}^{x}\rangle\dif u\right|\right],\\
&\mathcal{I}_3=\E\left[\left|\int_0^1(1-u)\nabla^2 f(u\tl{Y}_{\gamma}^{x}+(1-u)X_{\gamma}^{x})(\tl{Y}_{\gamma}^{x}-X_{\gamma}^{x})^{\otimes 2}\dif u\right|\right].
\end{align*}

It follows from Lemma \ref{one} that
\begin{align*}
\mathcal{I}_1\leq&\|\nabla f\|\|X_{\gamma}^{x}-\tilde{Y}_{\gamma}^{x}\|_2\leq C (1+|x|)\|\nabla f\|\gamma^{\frac{3}{2}},\\
\mathcal{I}_2\leq&\|\nabla^2 f\|\|X_{\gamma}^{x}-x\|_2\|X_{\gamma}^{x}-\tilde{Y}_{\gamma}^{x}\|_2\leq C(1+|x|^2)\|\nabla^2 f\|\gamma^{2},\\
\mathcal{I}_3\leq&C\|\nabla^2 f\|\|X_{\gamma}^{x}-\tilde{Y}_{\gamma}^{x}\|_2^2\leq C(1+|x|^2)\|\nabla^2 f\|\gamma^{3}.
\end{align*}
Collecting the above inequalities completes the proof.
\end{proof}

\begin{proof}[Proof of Lemma \ref{Step}]
    The proof of this lemma follows the same lines as the proof of \cite[Lemma 3.7]{jtp2} and is omitted.
\end{proof}

 \section*{Acknowledgements}
 The authors were partly supported by NSF of China with grants Nos: 12101392, 12271475 and U23A2064.

\end{document}